\documentclass[10pt]{article}
\usepackage[all]{xy}
\usepackage{enumerate, amsfonts, amsthm}

\newcommand{\mR}{\mathbb{R}}
\newcommand{\mC}{\mathbb{C}}
\newcommand{\mQ}{\mathbb{Q}}
\newcommand{\mZ}{\mathbb{Z}}
\newcommand{\mK}{\mathbb{K}}

\newtheorem{Def}{Definition}[section]
\newtheorem{Thm}[Def]{Theorem}

\newtheorem{Lem}[Def]{Lemma}

\begin{document}

\title{Twisted Equivariant {K}-theory for Proper actions of Discrete Groups}
\author{Christopher Dwyer}

\maketitle

\begin{abstract}
We will make a construction of twisted equivairant {K}-theory for proper actions of discrete groups by using ideas of L\"{u}ck and Oliver \cite{oliver-completion-theorem-2001} to expand a construction of Adem and Ruan \cite{adem-ruan}.
\end{abstract}

\section{Introduction}
Topological $K$-theory was first introduced in the 1960's as a generalized cohomology theory \cite{atiyah;k-theory}, and it was soon extended to the equivariant world by Segal \cite{segal;equivariant-k-theory}. In 1970, the notion of twisted $K$-theory was first introduced by Donovan and Karoubi \cite{donovan-karoubi;local-coefficients}. They twisted the $K$-theory of a space $X$ by a torsion element of $H^3(X,\mZ)$, and it wasn't until the 1980's that Rosenberg \cite{rosenberg;continuous-trace-algebras-bundle-theoretic} showed how to twist by a non-torsion element of $H^3(X,\mZ)$. He used the connections between $K$-theory and operator algebras to define the twisting.

In the 1990's, twisted $K$-theory gained notice when it came to the attention of physicists studying string theory. They found that certain charges can be described using the twisted $K$-theory of spacetime. More recently Freed, Hopkins, and Teleman proved that the twisted equivariant $K$-theory of a compact simple simply connected Lie group $G$ is the Verlinde algebra of the loop group of $G$ \cite{freed;twisted-verlinde-algebra,fht1,fht2}. 

Atiyah and Segal \cite{atiyah-segal;twisted} have recently written a paper which collects the definitions and constructions for twisted $K$-theory and twisted equivariant $K$-theory that are used in \cite{fht1,fht2} and \cite{bunke-schick;T-duality}. For the construction of twisted equivariant $K$-theory, they require that the group $G$ be compact. Motivated by orbifolds, Adem and Ruan \cite{adem-ruan} constructed twisted equivariant $K$-theory in a different manner for a finite group $G$. More general constructions can describe both notions, and work has been done by Jean-Louis Tu and Ping Xu \cite{tu-xu;chern-orbifolds} and also in their collaboration with Camille Laurent-Gengoux \cite{tu-xu-lg;differentiable-stacks} on this.

In this paper we will extend the construction of Adem and Ruan by following the example of L\"{u}ck and Oliver \cite{oliver-completion-theorem-2001} to define twisted equivariant $K$-theory for proper actions of general discrete groups. For simplicity, this paper is meant to be self contained, but in the interest of brevity, many of the elementary statements are left unproven. The beginning of every section contains references for a deeper treatment of the elementary materials. The paper is organized as follows:

Section 2 introduces projective representations of finite groups and some of the basic results involving them. Particular attention is given to the similarity between the theory of $\alpha$-twisted representations and complex representations. The subgroup theorem is particularly important as it will play a role in section 4. 

In section 3, twisted equivariant $K$-theory is constructed in the case of a torsion cocycle $\alpha$. The first part of the section is devoted to background and Theorem \ref{bundle}, which is a key component of the construction. The main theorem (Theorem \ref{main}) is stated and proved in the second part.

Section 4 is devoted to the definition of twisted Bredon cohomology and the description of a spectral sequence relating it to twisted equivariant $K$-theory. First, we define the coefficient functor along with a twisted product which is analogous to the one defined in section 3. Some properties of this twisted product are proved in section 3. Then, we will sketch the construction of the spectral sequence and use a result of L\"{u}ck \cite{luck} to show that it collapses at the $E_2$ page after tensoring with $\mC$.  A consequence is that this spectral sequence is a module over the non-twisted equivariant Atiyah-Hirzebruch spectral sequence in a natural way.

Section 5 talks about how this fits into the general theory of twisted equivariant $K$-theory coming from \cite{atiyah-segal;twisted}, \cite{tu-xu-lg;differentiable-stacks}, \cite{bunke-schick;T-duality}, and \cite{tu-xu;chern-orbifolds}. I owe great thanks to many people for helping me understand this material including but not limited to Alejandro Adem, Yongbin Ruan, my father, Wolfgang L\"{u}ck, Bob Oliver, Ian Leary, and Johann Leida. I also owe great thanks to the referee for suggesting a more elegant formulation of the results.

\section{Projective Representations}
This section will introduce some of the basic properties of projective representations. In particular, we will focus on the similarities of the theory of projective representations of a finite group with the theory of complex representations of a finite group. For this section, let $G$ be a finite group. More information about projective representations and proofs of the statements can be found in \cite{karpilovsky;projective-representations-finite-groups}.

\subsection{$\alpha$-twisted representations}

A projective representation of $G$ is a map which is not quite a complex representation but is only off by multiplication by a complex number.  
\begin{Def}
Let $V$ be a complex vector space. A \emph{projective representation} of $G$ is a map $\rho: G \rightarrow GL(V)$  such that there exists $\alpha: G \times G \rightarrow \mC^*$ with 
\begin{itemize}
\item[(i)] $\rho(x)\rho(y) = \alpha(x,y)\rho(xy)$ for all $x,y \in G$, and
\item[(ii)] $\rho(1)=I$.
\end{itemize}
\end{Def}

It is not hard to see that $\alpha$ is going to be a 2-cocycle of $G$, i.e. $\alpha \in Z^2(G, \mC^*)$. The associativity of multiplication in $GL(V)$ combined with condition $(i)$ gives the cocycle condition, and $\alpha(g,1) = \alpha(1,g) = 1$ for all $g \in G$ by condition $(ii)$. 

From now on we will focus on projective representations which are associated to the same cocycle $\alpha$. Any projective representation associated to $\alpha$ is called an \emph{$\alpha$-twisted} representation of $G$. In the case of complex representations of $G$, we have a notion of when two representations are isomorphic. We have a similar idea for $\alpha$-twisted representations. 

\begin{Def}
Two $\alpha$-twisted representations $\rho_i:G \rightarrow GL(V_i)$, $i=1,2$, are called 
linearly equivalent if there is a vector space isomorphism $f:V_1 \rightarrow V_2$ with
$$\rho_2(g) = f \rho_1(g) f^{-1} \, \, \, \mbox{ for all } g \in G$$
\end{Def}

Clearly the direct sum of two $\alpha$-twisted representations is again an $\alpha$-twisted representation. So we can form the monoid of linear isomorphism classes of $\alpha$-twisted representations of $G$, denoted by $M_\alpha(G)$. Let $R_\alpha(G)$ be its associated Grothendieck group, the $\alpha$-twisted representation group. Note that if $\alpha$ is the trivial cocycle, then $R_\alpha(G)$ is $R(G)$, the complex representation ring of $G$.

It is not hard to see that the tensor product of two $\alpha$-twisted representations is not an $\alpha$-twisted representation, but that doesn't mean that tensor products are uninteresting in this context. If $\alpha$ and $\beta$ are two cocycles, then the tensor product of an $\alpha$-twisted representation $\rho:G \rightarrow GL(V)$ and a $\beta$-twisted representation $\tau:G \rightarrow GL(W)$ is an $\alpha+\beta$-twisted representation. This can easily be seen by forming $\rho \otimes \tau: G \rightarrow GL(V \otimes W)$. Note that this can be extended to a pairing 
$$R_\alpha(G) \otimes R_\beta(G) \rightarrow R_{\alpha + \beta}(G).$$

In order to study $R_\alpha(G)$, we introduce the $\alpha$-\emph{twisted group algebra} $\mC^{\alpha}G$. We denote by $\mC^{\alpha}G$ the vector space over $\mC$ with basis \{$\bar{g}\}_{g \in G}$ with product $\bar{g} \cdot \bar{h} = \alpha(g, h)\overline{gh}$ extended by linearity. This makes $\mC^{\alpha}G$ into a $\mC$-algebra with $\bar{1}$ as the identity.

\begin{Def}
If $\alpha$ and $\beta$ are cocycles of $G$, then $\mC^{\alpha}G$ and $\mC^{\beta}G$ are \emph{equivalent} if there exists a $\mC$-algebra isomorphism $\phi:\mC^{\alpha}G \rightarrow \mC^{\beta}G$ and a mapping $t:G \rightarrow \mC^*$ such that $\phi(\bar{g}) = t(g)\tilde{g}$ for all $g \in G$, where {$\bar{g}$} and {$\tilde{g}$} are the bases for the two twisted group algebras.
\end{Def}

This defines an equivalence relation on such twisted algebras, and we have the following classification result 

\begin{Lem}
There is an equivalence of twisted group algebras, $\mC^{\alpha}G \simeq \mC^{\beta}G$, if and only if $\alpha$ is cohomologous to $\beta$. In fact, $\alpha \longmapsto \mC^{\alpha}G$ induces a bijective correspondence between $H^2(G, \mC^*)$ and the set of equivalence classes of twisted group algebras of $G$ over $\mC$.
\end{Lem}

In \cite[Theorem 2.2.5]{karpilovsky;projective-representations-finite-groups}, it is shown that these twisted group algebras play the same role in determining $R_\alpha(G)$ that $\mC G$ plays in determining $R(G)$.

\begin{Thm}
There is a bijective correspondence between $\alpha$-twisted representations of $G$ and $\mC^{\alpha}G$-modules. This correspondence preserves sums and bijectively maps linearly equivalent (irreducible, completely reducible) representations into isomorphic (irreducible, completely reducible) modules.
\end{Thm}

$\mC^{\alpha}G$ is also called the \emph{$\alpha$-twisted regular representation} of $G$. In \cite{karpilovsky;projective-representations-finite-groups}, it is shown that every $\mC^{\alpha} G$-module is projective.

\subsection{Induced $\alpha$-Representations}
Now let $\alpha$ be a fixed cocycle of $G$. We will use $\alpha$ to denote both the cocycle and its restriction to any subgroup as long as it is clear from the context which is meant. With this convention, for a subgroup $H$ of $G$ we can identify $\mC^{\alpha}H$ with the subalgebra of $\mC^{\alpha}G$ consisting of all $\mC$ linear combinations of the elements \{$\bar{h}|h \in H$\}. If $V$ is a $\mC^{\alpha}G$-module, then we shall denote by $V_H$ the $\mC^{\alpha}H$-module obtained by the restriction of the algebra. So $V_H$ equals $V$, but only the action of $\mC^{\alpha}H$ is defined on $V_H$. This process is called \emph{restriction} and it takes any $\mC^{\alpha}G$-module $V$ to a uniquely determined $\mC^{\alpha}H$-module $V_H$.

There is a dual process of \emph{induction}. Let $W$ be a $\mC^{\alpha}H$-module. We are considering $\mC^{\alpha}H$ as a subalgebra of $\mC^{\alpha}G$, and so we can define a $\mC^{\alpha}G$-module structure on the tensor product $\mC^{\alpha}G \otimes_{\mC^{\alpha}H} W$. We will denote this induced module $W^G$.

Here are some formal properties of induced modules

\begin{Lem}
Let $V_1$ and $V_2$ be submodules of a $\mC^{\alpha}H$-module $V$.
\begin{itemize}
\item[(i)] $V^G_1 \subseteq V^G_2$ if and only if $V_1 \subseteq V_2$.
\item[(ii)] $V^G_1 = V^G_2$ if and only if $V_1 = V_2$.
\item[(iii)] $(V_1 + V_2)^G = V^G_1 + V^G_2$.
\item[(iv)] $(V_1 \bigcap V_2)^G = V^G_1 \bigcap V^G_2$.
\item[(v)] If $V = V_1 \oplus V_2$, then $V^G = V^G_1 \oplus V^G_2$.
\end{itemize}
\end{Lem}

\begin{Lem}
Let $H$ be a subgroup of $G$ and let $\xymatrix@1{U \ar[r]^-{\lambda} & V \ar[r]^-{\mu} & W}$ be a sequence of homomorphisms of $\mC^{\alpha}H$-modules. Then the following properties hold:
\begin{itemize}
\item[(i)] The sequence
$$\xymatrix@1{
0 \ar[r] &  U \ar[r]^-{\lambda} & V \ar[r]^-{\mu} & W \ar[r] &  0 & & (*)
}$$
is exact if and only if the corresponding sequence of $\mC^{\alpha}G$-modules
$$\xymatrix@1{
0 \ar[r] & U^G \ar[r]^-{1 \otimes \lambda} &  V^G \ar[r]^-{1 \otimes \mu} & W^G \ar[r] &  0  & & (**)
}$$
is exact.
\item[(ii)] Suppose that the sequence \emph{(*)} is exact. Then \emph{(*)} splits if and only if \emph{(**)} splits.
\item[(iii)] If $U$ is a submodule of $V$, then $V^G/U^G \cong (V/U)^G$.
\end{itemize}
\end{Lem}

Some notation is needed for the following theorem. Let $H$ be a subgroup of $G$, and let $V$ be a $\mC^{\alpha}H$-module. We can always choose a $\mC^{\alpha}G$-module $W$ so that $V$ is a submodule of $W_H$. Then, for any $g \in G$, $\bar{g}V = \{\bar{g}v|v \in V \}$ is a $\mC$-subspace of $W$. In fact, $\bar{g}V$ is a $\mC^{\alpha}(gHg^{-1})$-module  whose isomorphism class is independent of the choice of $W$. Now define $V^{(g)}$ to be the $\mC^{\alpha}(gHg^{-1})$-module whose underlying space is $V$ and with $\bar{x}$, for $x \in gHg^{-1}$, acting as 
$$\begin{array}{rcl}
\bar{x} \cdot v & = & (\bar{g}^{-1} \bar{x} \bar{g})v \\
    & = &\alpha^{-1}(xg, g^{-1})\alpha^{-1}(g, g^{-1}xg)\alpha(g, g^{-1})\overline{g^{-1}xg}v. \\
\end{array}$$
Note that if $g \in C_G(H)$, then $V^{(g)}$ is isomorphic to $V$ as $\mC^{\alpha}H$-modules since $\overline{g^{-1}xg} = \bar{x}$.

\begin{Thm}[Subgroup Theorem] \label{sub}
Let $H$ and $K$ be subgroups of the finite group $G$. Let $T$ be a set of double coset representatives for $(K, H)$ in $G$, and for $\alpha \in Z^2(G, \mC^*)$, let $V$ be a $\mC^{\alpha}H$-module. For each $t \in T$, denote by $V_t$ the restriction of $V^{(t)}$ to $\mC^{\alpha}(tHt^{-1} \bigcap K)$. Then
$$(V^G)_K \cong \bigoplus_{t \in T} (V_t)^K$$
as $\mC^{\alpha}K$-modules.
\end{Thm}

\section{Twisted Equivariant $K$-theory}
We will construct twisted equivariant $K$-theory out of finite rank equivariant bundles similar to \cite{segal;equivariant-k-theory}. However, we only consider equivariant bundles whose underlying structure is a complex vector bundle over $X$. Another concern is that in \cite{segal;equivariant-k-theory} the group must be compact. Since we are interested in actions of general discrete groups, we can no longer assume this. However, we will assume that the group is acting properly. This means that all the subgroups that fix a point are compact and therefore finite. We will using the work of L\"{u}ck and Oliver \cite{oliver-completion-theorem-2001} to construct the \emph{$\alpha$-twisted} equivariant $K$-theory of a finite proper $G$-CW pair for a discrete group $G$.

First we need some background for the construction. Let $\alpha \in Z^2(G, U(1))$ be a torsion cocycle of order $n$. Since 
$\alpha$ is torsion of order $n$, $\alpha^n$ is cohomologous to the 
trivial cocycle, i.e., there exists a cocycle $t \in Z^1(G, U(1))$ with 
$\alpha^n = (\delta t)$. Define a cocycle $u \in Z^1(G, U(1))$ by 
$u(g) = (t(g))^{-1/n}$. We can define such a  $u$ since $U(1)$ is a divisible group. Using multiplicative notation, this means we can always take $n^{th}$ roots. $\beta = \alpha (\delta u)$ is again torsion of 
order $n$,  and $\beta$ takes values in $\mZ/n\mZ$. So by restricting to a 
cohomologous cocycle if necessary, we may assume that $\alpha$ takes 
values in $\mZ/n\mZ$. We call such a cocycle \emph{normalized}.

So, we can take $\alpha$ to represent a central extension of $G$ by $\mZ/n\mZ$ as below.
$$1 \longrightarrow \mZ/n\mZ \longrightarrow G_{\alpha} \stackrel{\rho}{\longrightarrow} G \longrightarrow 1$$

Recall that elements of $G_{\alpha}$ are of the form $(g, \sigma^k)$ for $g \in G$, $\sigma$ the generator of $\mZ/n\mZ$ and $0 \leq k \leq n-1$; where multiplication is given by $(g, \sigma^j) \cdot (h, \sigma^k) = (gh, \alpha(g,h) \sigma^{j+k})$. We can form $\mC G_{\alpha}$, the algebra over $\mC$ with basis  $\{ \overline{(g, \sigma^j)} | g \in G, \, 0 \leq j \leq n-1 \}$. We also have $\mC^{\alpha} G$, the $\alpha$-twisted group algebra of $G$, and  $\mC^{\alpha} G$ has basis $\{ \tilde{g} | g \in G \}$ as an algebra over $\mC$. 

For the following let $\zeta = e^{2 \pi i/n}$. There is a natural $\mC$-algebra monomorphism $\phi:\mC^{\alpha}G \rightarrow \mC G_{\alpha}$ with
$$\phi(\tilde{g}) = \displaystyle\frac{1}{n} \displaystyle\sum_{j=0}^{n-1} \zeta^{n-j} \overline{(g, \sigma^j)}$$
extended by linearity. This is a $\mC$-algebra homomorphism because if $\alpha(g, h) = \sigma^j$, then
$$\begin{array}{rcl}
\phi(\tilde{g})\phi(\tilde{h}) & = & \displaystyle\frac{1}{n^2}\left(\displaystyle\sum_{0 \leq k \leq n-1} \zeta^{n-k}\overline{(g, \sigma^k)}\right) \left(\displaystyle\sum_{0 \leq m \leq n-1} \zeta^{n-m}\overline{(h, \sigma^m)}\right) \\
 & & \\
 & = & \displaystyle\frac{1}{n^2}\left(n \overline{(gh, \sigma^j)} + n\zeta^{n-1}\overline{(gh, \sigma^{j+1})} + \cdots + n\zeta^1\overline{(gh, \sigma^{j-1})}\right)\\
 & & \\
 & = &  \alpha(g,h)\phi(\widetilde{gh}) \\
\end{array}$$ 
Note that $\phi(\zeta^{k} v) = \overline{(1, 
\sigma^k)}\cdot\phi(v)$ for any $v \in \mC^{\alpha} G$. Since $\phi$ 
takes the basis elements of $\mC^\alpha G$ to linearly independent 
vectors in $\mC G_{\alpha}$, $\phi$ is clearly a monomorphism. This 
gives a way to identify $\mC^{\alpha} G$ with a subspace of $\mC 
G_{\alpha}$. This identification is natural with respect to restriction 
to subgroups as follows. For any group monomorphism $f:H \rightarrow G$, we have a commutative diagram 
$$\begin{array}{rcl}
\mC^{\alpha} H  & \stackrel{\tilde{f}}{\longrightarrow} & \mC^{\alpha} G \\
\downarrow \phi & & \phi \downarrow \\
\mC H_{\alpha} & \stackrel{\bar{f}}{\longrightarrow} & \mC G_{\alpha} \\
\end{array}$$

It is well known that if $[G_{\alpha}:H_{\alpha}] = k$, then $\mC G_{\alpha} \cong 
\oplus_{i=1}^{k} \mC H_{\alpha}$ as right $\mC H_\alpha$-modules by the 
relationship 
$$
\mC G_{\alpha} = \overline {t_1} (\mC H_{\alpha}) \oplus \overline{t_2} 
(\mC H_{\alpha}) \oplus \overline{t_3} (\mC H_{\alpha}) \oplus \cdots 
\oplus \overline{t_k} (\mC H_{\alpha})
$$
where $\{ t_1, t_2, t_3, \ldots, t_k \}$ is a left transversal for 
$H_{\alpha}$ in $G_{\alpha}$. We may assume that $t_1$ is the 
identity element of $G_{\alpha}$.  Note that $\{ \rho(t_1), \rho(t_2), 
\rho(t_3), \ldots, \rho(t_k) \}$ is a left transversal for $H$ in $G$, where $\rho:G_\alpha \rightarrow G$ is the projection. So we have
$$ \mC^{\alpha}G = \bigoplus_{i=1}^{k} \widetilde{\rho(t_i)}(\mC^{\alpha}H)$$

This fact combined with the canonical inclusion above gives the 
following commutative diagram of $H_{\alpha}$-modules
$$\begin{array}{rcl}
\bigoplus_{i=1}^{k} \widetilde{\rho(t_i)}(\mC^{\alpha}H) & \stackrel{=}{\longrightarrow} & 
\mC^{\alpha}G  \\
\downarrow \oplus \phi & & \downarrow \phi  \\
\bigoplus_{i=1}^k \overline{t_i}(\mC H_{\alpha}) & \stackrel{=}{\longrightarrow} & \mC 
G_{\alpha} \\
\end{array}$$
This will be important in the proof of the following theorem. 

\begin{Thm} \label{bundle}
Let $G$ be a discrete group, and let $\alpha$ be a normalized torsion cocycle of $G$ as above. Then for any finite proper $G$-CW complex $X$, there is a $G_{\alpha}$-bundle $E$ over $X$ with fiber $E|_x \cong (V_{\alpha, G_x})^{\oplus k}$ for some $k$ where $V_{\alpha, G_x} = \mC^{\alpha}G_x$. That is, the fiber of $E$ is the direct sum of a number of copies of the twisted regular representation of $G_x$.
\end{Thm}
\noindent \emph{Proof.}
First, we need to note that since $G_{\alpha}$ is a finite extension of $G$, we can think of $X$ as a finite, proper $G_\alpha$-CW complex as follows. For any subgroup $H_i$ of $G$, $\alpha$ restricts to the subgroup $H_i$ to define an extension $H_{i, \alpha}$ of $H_i$ by $\mZ/n\mZ$. Then for every $G$-equivariant $n$-cell of $X$, say $(G/H_i \times D^n, \, G/H_i \times S^{n-1})$, define a $G_\alpha$-equivariant $n$-cell of $X$, $(G_{\alpha}/H_{i, \alpha} \times D^n, \, G_{\alpha}/H_{i, \alpha} \times S^{n-1})$. Gluing these cells together in the obvious way gives a $G_\alpha$-CW structure for $X$.

Another fact we will need is that there is a $G_\alpha$- bundle $E^\prime$ over $X$ with fiber $E'|_x$ the direct sum of a number of copies of the regular representation of $G_{x,\alpha}$ \cite[Corollary 2.7]{oliver-completion-theorem-2001}. The bundle we want will be a subbundle of $E^\prime$.

Using the notation above, we will write elements of $G_\alpha$ in the form $(g, \sigma^i)$ with $g \in G$ and $\sigma$ the generator of $\mZ/n\mZ$, and let $\zeta = e^{2\pi i /n}$. Note that since $(1,\sigma^j)$ fixes $X$ the action of $(1,\sigma^j)$ on $E^\prime$ must take any fiber to itself. Define $E$ to be the subbundle of $E^\prime$ consisting of all vectors $v$ with $(1,\sigma^j) \cdot v = \zeta^j v$ for $0 \leq j \leq n-1$.  

It still remains to be shown that $E$ is the bundle we want. The fiber $E|_x$ is going be a subspace of the fiber $E^\prime |_x = (\mC G_{x,\alpha})^{\oplus k}$, so $E|_x$ is going to contain $(\mC^\alpha G_x)^{\oplus k}$ in an obvious way. We need to show that that there isn't anything else. So let $v$ be a vector in $E|_x$. We can view $v$ as a vector in $E^\prime |_x = (\mC G_{x,\alpha})^{\oplus k}$, and we can express $v$ in the following way,
$$v = \displaystyle\sum_{j=1}^{k} \displaystyle\sum_{h^j \in G_x^j} \displaystyle\sum_{l =0}^{n-1} a_{h,l}^j\overline{(h^j, \sigma^l)}
$$
where $G_x^j$ is the $j^{th}$ copy of $G_x$ in the fiber. A simple calculation shows that since $(1,\sigma^j) \cdot v = \zeta^j v$, $v \in (\mC^\alpha G_x)^{\oplus k}$. So $E$ is the bundle we want.
\hfill \qed

\vspace{0.2cm}

Let $G$ be a discrete group and let $\alpha \in Z^2(G, U(1))$ be a normalized
cocycle of finite order $n$, as above. We know that $\alpha$ defines a central 
extension of $G$ by $U(1)$, and this extension factors through a central extension of $G$ by $\mZ/n\mZ$ in the following sense. In the diagram below, $G_\alpha$ is included as a subgroup of $\widetilde{G_\alpha}$, and the diagram commutes.

$$\begin{array}{rcccccccl}
1 & \rightarrow & \mZ/n\mZ & \rightarrow & G_\alpha & \rightarrow & G & \rightarrow & 1 \\
\left| \right| & & \downarrow & & \downarrow & & \left| \right| & & \left| \right| \\
1 & \rightarrow & U(1) & \rightarrow & \widetilde{G_\alpha} & \rightarrow & G & \rightarrow & 1 \\
\end{array}$$

From the above diagram,  it is not hard to see that any $\widetilde{G_{\alpha}}$-equivariant bundle in which the central $U(1)$ acts by complex multiplication on the fibers restricts to a $G_{\alpha}$-equivariant bundle in which the central $\mZ/n\mZ$ acts by multiplication by powers of $e^{2\pi i/n}$ on the fibers. Equivalently, it is clear that for any $G_{\alpha}$-bundle in which the central $\mZ/n\mZ$ acts by complex multiplication on the fibers, the action can be extended to a $\widetilde{G_{\alpha}}$ action with the central $U(1)$ acting by complex multiplication on the fibers. So it suffices to consider $G_{\alpha}$-equivariant bundles with this action.

\begin{Def} Let $G$ be a discrete group, $X$ a finite, proper $G$-CW 
complex, and $\alpha \in Z^2(G, U(1))$ a normalized torsion cocycle of 
order $n$. An \emph{$\alpha$-twisted $G$-bundle} over $X$ is a complex vector 
bundle $E \rightarrow X$ with an action of $\widetilde{G_\alpha}$ on $E$ which covers the action of $G$ on $X$ and which restricts to the complex multiplication action of $U(1) \subset \widetilde{G_\alpha}$ on the fibers.
\end{Def}

So far, we have only been considering \emph{torsion} cocycles. It turns out that for our purposes, it suffices to consider only torsion cocycles, as the following lemma will show.

\begin{Lem}
Let $X$ and $G$ be as in the definition, and let $E$ be an $\alpha$-twisted $G$ bundle of rank $n$. Then the order of $\alpha$ divides $n$.
\end{Lem}
\noindent \emph{Proof.}
Let $E$ be an $\alpha$-twisted $G$ bundle of rank $n$ over $X$. Since $E$ is a $(\widetilde{G_\alpha}, GL_n(\mC))$ bundle, we will get a homomorphism $\tilde{f}:\widetilde{G_\alpha} \rightarrow GL_n(\mC)$ by restricting to an orbit and trivializing the bundle there. $\tilde{f}$ restricts to a homomorphism $f:G \rightarrow PGL_n(\mC)$ after dividing out by $U(1)$. The map $f$ represents an element of $H^1(G; PGL_n(\mC))$ in an obvious way. 

There is a long exact sequence in group cohomology coming from the short exact sequence of coefficients
$$\begin{array}{rcccccccl}
1 & \rightarrow & U(1) & \rightarrow & GL_n(\mC) & \rightarrow & PGL_n(\mC) & \rightarrow & 1. 
\end{array}$$
This gives a map $\delta: H^1(G;PGL_n(\mC)) \rightarrow H^2(G; U(1))$ which, it is easy to see, takes the cohomology class of $f$ to the cohomology class of $\alpha$. Because of the following commuting diagram of coefficients
$$\begin{array}{rcccccccl}
1 & \rightarrow & \mZ/n\mZ & \rightarrow & SL_n(\mC) & \rightarrow & PGL_n(\mC) & \rightarrow & 1 \\
& & \downarrow & & \downarrow & & \left| \right| & & \\
1 & \rightarrow & U(1) & \rightarrow & GL_n(\mC) & \rightarrow & PGL_n(\mC) & \rightarrow & 1, \\
\end{array}$$
we have the following commutative diagram in cohomology
$$\xymatrix{
H^1(G; PGL_n(\mC)) \ar[r]^{\delta^\prime}  \ar[d]_{\cong} & H^2(G; \mZ/n\mZ) \ar[d] \\
H^1(G; PGL_n(\mC)) \ar[r]^\delta & H^2(G; U(1)). }
$$
Therefore, $\alpha$ lies in the image of $H^2(G; \mZ/n\mZ) \rightarrow H^2(G; U(1))$ and so it must be $n$-torsion. This proves the statement.\hfill \qed

\vspace{0.2cm}

In other words, if $\alpha$ is not a torsion cocycle, there are no finite rank $\alpha$-twisted $G$ bundles. If $\alpha$ is torsion, then we know there is at least one $\alpha$-twisted bundle by Theorem \ref{bundle}. 

Now if $E$ and $F$ are both $\alpha$-twisted $G$-bundles, then so is $E \oplus F$. This means that isomorphism classes of $\alpha$-twisted bundles over $X$ form a monoid. Let $\,^{\alpha}K_G(X)$ be the associated Grothendieck group. It is important to note that if $E$ and $F$ are both $\alpha$-twisted $G$-bundles, $E \otimes F$ is {\bf not} an $\alpha$-twisted $G$-bundle. So $\,^{\alpha}K_G(X)$ doesn't have a ring structure.

The following observation will allow us to extend the construction of $\,^\alpha K_G(X)$ to a $\mZ/2\mZ$-graded equivariant cohomology theory. Any $G_\alpha$-equivariant bundle $E$ over $X$ will determine a representation of $\mZ/n\mZ$ for each $G$-connected component of $X$. This comes from the fact that $\mZ/n\mZ$ acts trivially on $X$. So each fiber of $E$ gives a representation of $\mZ/n\mZ$, and because $E$ is a bundle, the representation is fixed over every $G$-connected component of $X$. Let $\mK_{G_\alpha}(X)$ be the Grothendieck group of the monoid of isomorphism classes of $G_\alpha$-equivariant bundles over $X$. Then, if $X$ is $G$-connected, we have a splitting
\begin{equation}
\label{splitting} \mK_{G_\alpha}(X) =  \bigoplus_{V \in Irr(\mZ/n\mZ)} \mK_{G_\alpha,V}(X). 
\end{equation}
where $Irr(\mZ/n\mZ)$ is the set of isomorphism classes of irreducible complex representations of $\mZ/n\mZ$. Here $\mK_{G_\alpha,V}(X)$ consists of the $G_\alpha$-equivariant vector bundles which restrict to the representation $V$ in each fiber. The splitting comes from the fact that every complex representation of $\mZ/n\mZ$ splits into a direct sum of irreducibles. 

So we have a splitting according to the irreducible complex characters of $\mZ/n\mZ$,  and it is clear that $\,^\alpha K_G(X)$ is the direct summand corresponding to the representation given by taking the generator of $\mZ/n\mZ$ to multiplication by $e^{2 \pi i /n}$. Also, the splitting is functorial in the category of finite, proper $G$-CW complexes and $G$-maps via pullback. Given any $G$-map $f:X \longrightarrow Y$ between $G$-CW complexes $X$ and $Y$, $f$ is trivially $G_\alpha$-equivariant, and the induced map $f^*:\mK_{G_\alpha}(Y) \longrightarrow \mK_{G_\alpha}(X)$ will respect the splitting (\ref{splitting}). So $f^*$ will take $\mK_{G_\alpha, V}(Y)$ to $\mK_{G_\alpha,V}(X)$ for each $V \in Irr(\mZ/n\mZ)$. Also, as noted for $\,^\alpha K_G(X)$, the Whitney sum of bundles will restrict to an operation on each summand in the splitting (\ref{splitting}). So each summand is an abelian group with respect to Whitney sum. Tensor product of bundles will not restrict to an operation on each summand, however.

Since $G_\alpha$ is a finite extension of $G$, it is again a discrete group, and $X$ is a proper $G_\alpha$-CW complex if it is a proper $G$-CW complex as seen in Theorem \ref{bundle}. Therefore we can describe the $G_\alpha$-equivariant $K$-theory of $X$ using the results of L\"{u}ck and Oliver \cite{oliver-completion-theorem-2001}. In particular, if we again let $\mK_{G_{\alpha}}(X)$ denote the Grothendieck group of the monoid of isomorphism classes of $G_{\alpha}$-equivariant vector bundles on $X$, then we define the $G_\alpha$-equivariant $K$ groups of $X$ by
$$K_{G_\alpha}^{-n}(X) = Ker[\mK_{G_\alpha}(X \times S^n) \stackrel{incl^*}{\longrightarrow} \mK_{G_\alpha}(X)].$$
For a finite, proper $G$-CW pair $(X,A)$, define 
$$K_{G_\alpha}^{-n}(X, A) = Ker[\mK_{G_\alpha}^{-n}(X \cup_A X) \stackrel{i^*_2}{\longrightarrow} \mK_{G_\alpha}^{-n}(X)].$$
L\"{u}ck and Oliver \cite{oliver-completion-theorem-2001} show that these definitions give a $\mZ/2\mZ$ graded equivariant cohomology theory and Bott periodicity holds in this context.

We can now state the main theorem of the paper. 

\begin{Thm}[Main Theorem] \label{main}
 Given a discrete group $G$ and a normalized torsion cocycle $\alpha \in Z^2(G, U(1))$, the groups $\,^{\alpha}K_G^{-n}(X,A)$ extend to a $\mZ/2\mZ$-graded equivariant cohomology theory called $\alpha$-twisted equivariant $K$-theory on the category of finite proper $G$-CW pairs. If $\alpha$ is the trivial cocycle, then $\,^{\alpha}K^*_G(X,A)$ is just the equivariant $K$-theory defined by L\"{u}ck and Oliver \cite{oliver-completion-theorem-2001}. For two normalized torsion cocycles $\alpha$ and $\beta$, there is a twisted product 
$$\,^{\alpha}K_G^*(X, A) \bigotimes \,^{\beta}K_G^*(X, A) \rightarrow \,^{\alpha + \beta}K_G^*(X, A).$$
which gives $\,^{\alpha}K_G^*(X,A)$ a natural graded $K^*_G(X,A)$-module structure. Also, for finite $H \subseteq G$, there are natural isomorphisms $\,^{\alpha}K_G^0(G/H) \cong R_{res^G_H \alpha}(H)$ and $\,^{\alpha}K_G^1(G/H) = 0$. If $G$ is finite, then this construction agrees with the construction of Adem and Ruan \cite{adem-ruan}. 
\end{Thm}

The proof of Theorem \ref{main} will be the focus of the rest of this
section. The main issue will be proving the existence of the twisted product. The fact that we have an equivariant cohomology theory will follow from the fact that $\,^{\alpha}K_G(X)$ is a direct summand of $K_{G_{\alpha}}(X)$ and all the maps in \cite{oliver-completion-theorem-2001} respect the splitting in (\ref{splitting}). 

First, we see from the definition that each $K^{-n}_{G_\alpha}(X)$ splits as in (\ref{splitting}). Because the splitting is functorial in the category of finite, proper $G$-CW complexes and $G$-maps, we see that for any $G$-CW complex $X$, we can define $\,^\alpha K^{-n}_G(X)$ as the direct summand of $K^{-n}_{G_\alpha}(X)$ corresponding to the representation which takes the generator of $\mZ/n\mZ$ to multiplication by $e^{2 \pi i /n}$. With this definition, the extension of the groups $\,^{\alpha} K_G(X)$ to a $\mZ/2\mZ$-graded $G$-equivariant cohomology theory follows from the fact that the groups $K^{-n}_{G_\alpha}(X)$ extend to a $\mZ/2\mZ$-graded cohomology theory on the category of finite, proper $G$-CW pairs \cite{oliver-completion-theorem-2001}.

The properties of excision, homotopy invariance, disjoint union, and the long exact and Mayer-Vietoris sequences for $\,^{\alpha}K^{-n}_G(X)$ all follow directly from the same properties of $K^{-n}_{G_\alpha}(X)$ proven in \cite{oliver-completion-theorem-2001} by restricting all the maps to the appropriate direct summands. All of the maps in the constructions which are induced from $G_\alpha$-equivariant maps automatically respect the splitting (\ref{splitting}), and so restrict to the direct summands. The only map which needs to be checked is the connecting homomorphism in the Mayer-Vietoris sequence
$$d^{-n}:K^{-n}_{G_\alpha}(A) \rightarrow K^{-n+1}_{G_\alpha}(X).$$
However, it is easily seen to respect the splitting by its definition. It is the restriction of the induced map
$$\mK_{G_\alpha}((X_1 \times D^n) \bigcup_{A \times S^{n-1}} (X_2 \times D^n)) \stackrel{incl^*}{\longrightarrow} \mK_{G_\alpha}(X \times S^n)$$
to a subgroup of $\mK_{G_\alpha}((X_1 \times D^n) \bigcup_{A \times S^{n-1}} (X_2 \times D^n))$. Here $X_1$ and $X_2$  are part of a pushout diagram of finite proper $G$-CW complexes
$$\xymatrix{
A \ar[d]_{i_2} \ar[r]^{i_1} & X_1  \ar[d]^{j_1}\\
X_2 \ar[r]^{j_2} & X }
$$ 
where $i_1$ and $j_2$ are inclusions of subcomplexes.  

We should say a little about induction. Let $(X,A)$ be a finite proper $G$-CW pair, and let $H$ be a subgroup of $G$. From \cite[Lemma 3.4]{oliver-completion-theorem-2001}, we have 
$$i^{G_\alpha}_{H_\alpha}: K^{-n}_{H_\alpha}(X,A) \stackrel{\cong}{\longrightarrow} K^{-n}_{G_\alpha}(G_\alpha \times_{H_\alpha}(X,A)).$$
Here we are using the abuse of notation that $H_\alpha$ corresponds to the extension of $H$ by $\mZ/n\mZ$ which comes from restricting the cocycle $\alpha$ to the subgroup $H$. This isomorphism respects the splittings of $K^{-n}_{H_\alpha}(X,A)$  and $K^{-n}_{G_\alpha}(G_\alpha \times_{H_\alpha} (X,A))$. This follows from the construction of the isomorphism since given an $H_\alpha$-equivariant vector bundle $E$ and the corresponding $G_\alpha$-equivariant vector bundle $G_\alpha \times_{H_\alpha} E$, they will both clearly restrict to the same representation of $\mZ/n\mZ$ in the fibers over the points $x$ and $G_\alpha \times_{H_\alpha} x$. So the isomorphism will restrict to the appropriate summand, and induciton will hold for $\,^{\alpha}K^{-n}_G(X,A)$.

The Bott periodicity map will also respect the splitting (\ref{splitting}), since it involves the exernal tensor product with a bundle over $S^2$. So for any finite proper $G$-CW pair $(X,A)$ and normalized torsion cocycle $\alpha$ of $G$, we can define $\,^{\alpha} K^0_G(X,A)$ as the direct summand of $K^0_{G_\alpha}(X,A)$ corresponding to the representation of $\mZ/n\mZ$ which takes the generator to multiplication by $e^{2 \pi i/n}$. We can define $\,^\alpha K^1_G(X,A)$ as the corresponding direct summand of $K^1_{G_\alpha}(X,A)$. This gives a $\mZ/2\mZ$-graded  equivariant cohomology theory on the category of finite proper $G$-CW pairs which has the following  exact rectangle 
$$\xymatrix{
\,^{\alpha}K^0_G(X,A) \ar[r] & \,^{\alpha}K^0_G(X) \ar[r] &
\,^{\alpha}K^0_G(A) \ar[d]^{\delta^0}\\
\,^{\alpha}K^1_G(A) \ar[u]^{\delta^{-1}} & \ar[l] \,^{\alpha}K^1_G(X) & \ar[l] \,^{\alpha}K^1_G(X,A) }
$$
for any finite proper $G$-CW pair $(X,A)$. The rectangle again just comes from the rectangle we get from the $G_\alpha$-equivariant $K$-theory of the pair $(X,A)$ \cite[Theorem 3.2]{oliver-completion-theorem-2001}.
 
The only part of Theorem \ref{main} left to prove is the existence of an associative graded commutative twisted product. This will come from the graded commuative product constructed by L\"{u}ck and Oliver in \cite{oliver-completion-theorem-2001}. Given two normalized torsion cocycles $\alpha$ and $\beta$ with values in $\mZ/n\mZ$ and $\mZ/m\mZ$ respectively, we can form an exterior product $\alpha \oplus \beta$ which is a cocycle with values in $\mZ/n\mZ \times \mZ/m\mZ$. This defines an extension of $G$ by $\mZ/n\mZ \times \mZ/m\mZ$ which we can call $G_{\alpha \oplus \beta}$.

There are obvious projection maps $p_\alpha: G_{\alpha \oplus \beta} \rightarrow G_\alpha$  and $p_\beta: G_{\alpha \oplus \beta} \rightarrow G_\beta$, and these give rise to maps in equivariant $K$-theory
$$\phi_\alpha:K_{G_\alpha}(X) \rightarrow K_{G_{\alpha \oplus \beta}}(X), \,\,\, \phi_\beta:K_{G_\beta}(X) \rightarrow K_{G_{\alpha \oplus \beta}}(X)$$
for any $G$-space $X$. Note that if $E$ is a $G_{\alpha}$-bundle, then the central subgroup $1 \times \mZ/m\mZ \subset G_{\alpha \oplus \beta}$ acts trivially on $\phi_\alpha(E)$, and similarly if $F$ is a $G_\beta$-bundle, then the central subgroup $\mZ/n\mZ \times 1 \subset G_{\alpha \oplus \beta}$ acts trivially on $\phi_\beta(F)$.

We can form the composition 
$$\xymatrix@1{ 
K_{G_\alpha}(X) \bigotimes K_{G_\beta}(X) \ar[r]^-{\phi_\alpha \otimes \phi_\beta}  & K_{G_{\alpha \oplus \beta}}(X) \bigotimes  K_{G_{\alpha \oplus \beta}}(X) \ar[r]^-{\mu} & K_{G_{\alpha \oplus \beta}}(X) \,\,\,\,\,\,\,\,\,\,\,(3)
} \label{comp}$$
where $\mu$ is the product defined in \cite{oliver-completion-theorem-2001}. This will give us the associative graded commutative twisted product we want. Let us call this composition $\eta$.

Notice that $K_{G_{\alpha \oplus \beta}}(X)$ has a splitting just as in (\ref{splitting}). In this case, the direct summands will correspond to irreducible characters of $\mZ/n\mZ \times \mZ/m\mZ$. However, any irreducible character of $\mZ/n\mZ \times \mZ/m\mZ$ can be seen as the product of an irreducuble character of $\mZ/n\mZ$ with an irreducible character of $\mZ/m\mZ$. 

If we restrict $\eta$ to $\,^\alpha K_G(X) \bigotimes \,^\beta K_G(X)$, the image lands in the direct summand of $K_{G_{\alpha \oplus \beta}}(X)$ which corresponds to the representation of $\mZ/n\mZ \times \mZ/m\mZ$ which takes the generator of $\mZ/n\mZ \times 1$ to multiplication by $e^{2 \pi i / n}$ and the generator of  $1 \times \mZ/m\mZ$ to multiplication by $e^{2 \pi i/m}$. We will denote this image by $\,^{\alpha,\beta} K_{G_{\alpha \oplus \beta}}(X)$.  

\begin{Lem}
Let $\alpha, \, \beta \in Z^2(G, U(1))$ be two normalized torsion cocycles of the discrete group $G$. For any finite $G$-CW complex $X$, there is a canonical map
$$\varphi:\mK_{G_{\alpha + \beta}}(X) \longrightarrow \mK_{G_{\alpha \oplus \beta}}(X)$$ 
such that the restriction
$$\overline{\varphi}: \,^{\alpha + \beta} K_G(X) \longrightarrow \,^{\alpha,\beta} K_{G_{\alpha \oplus \beta}}(X)$$
is an isomorphism. Here $\alpha + \beta \in Z^2(G,U(1))$ represents the normalized torsion cocycle which is the sum of the cocycles $\alpha$ and $\beta$.
\end{Lem}
\noindent \emph{Proof.  }
Since they are normalized torsion cocycles, let $\alpha$ take values in $\mZ/n\mZ \subset U(1)$ and let $\beta$ take values in $\mZ/m\mZ \subset U(1)$.  Since $\alpha$ and $\beta$ are normalized, $\alpha + \beta$ will be normalized as well, and it will take values in $\mZ/k\mZ$ for some $k$ dividing the least common multiple of $n$ and $m$. In fact, $\mZ/k\mZ$ will be the image of $\mZ/n\mZ \times \mZ/m\mZ$ under the multiplication map
$$U(1) \times U(1) \rightarrow U(1).$$ 

There is an obvious projection map $p:G_{\alpha \oplus \beta} \rightarrow G_{\alpha + \beta}$ given by $(g,a,b) \mapsto (g,ab)$. Let $\varphi$ be the map induced by $p$ on equivariant $K$-theory.  Under this construction, if $E$ is an $(\alpha+\beta)$-twisted $G$-bundle, then $\varphi(E)$ will be the equivariant $G_{\alpha \oplus \beta}$-bundle with the same total space $E$, where $G_{\alpha \oplus \beta}$ acts on $\varphi(E)$ by the projection $p$ and the action of $G_{\alpha+\beta}$ on $E$. This action shows that $\varphi(E)$ is an element of $\,^{\alpha,\beta}K_{G_{\alpha \oplus \beta}}(X)$ from the definition, and it makes sense to talk about the restriction $\overline{\varphi}$. It is clear from this construction that at the level of bundles, $\overline{\varphi}$ will take two non-isomorphic $\alpha+\beta$-twisted bundles $E$ and $F$ to non-isomorphic $G_{\alpha \oplus \beta}$-equivariant bundles.

It remains to be shown that $\overline{\varphi}$ is a surjection. However, this is simple to see at the level of bundles. If we take a $G_{\alpha \oplus \beta}$-equivariant bundle $E$ with the generator of the central $\mZ/n\mZ \times 1$ (resp, central $1 \times \mZ/mZ$) acting by multiplication by $e^{2 \pi i/n}$ (resp. $e^{2 \pi i/m}$) on the fibers, then we can define a bundle $F$ with the same total space and an action of $G_{\alpha + \beta}$ as follows. The generator of the central $\mZ/k\mZ$ acts by multiplication by $e^{2 \pi i/k}$ on the fibers. Given an element $g \in G$, we will let $(g,1) \in G_{\alpha+\beta}$ act on $F$ the same way $(g,1,1) \in G_{\alpha \oplus \beta}$ acts on $E$. This action will be well defined, and it will make $F$ into an $\alpha + \beta$-twisted $G$-bundle. This construction means that every $G_{\alpha \oplus \beta}$-bundle with the given action corresponds to an $\alpha+\beta$-twisted $G$-bundle, and so $\overline{\varphi}$ is an isomorphism.
\hfill \qed

\vspace{0.2cm}

Using the above lemma, we can construct the product
$$\overline{\varphi}^{-1} \circ \eta: \,^{\alpha} K_G(X) \bigotimes \,^{\beta}K_G(X) \longrightarrow \,^{\alpha + \beta}K_G(X).$$
This product will extend to a product
$$\tau: \,^{\alpha}K^i_G(X) \bigotimes \,^{\beta}K^j_G(X) \longrightarrow \,^{\alpha + \beta}K^{i+j}_G(X)$$
in the same way that the product in \cite{oliver-completion-theorem-2001} extends. In fact, using the splitting (\ref{splitting}) 
$$K_{G_{\alpha \oplus \beta}}(X) = \bigoplus_{V \in Irred(\mZ/n\mZ \times \mZ/m\mZ)} K_{G_{\alpha \oplus \beta}, V}(X)$$
we have for the $G_{\alpha \oplus \beta}$-equivariant $K$-theory of a finite $G$-CW complex $X$, the product
$$K_{G_{\alpha \oplus \beta}}(X) \bigotimes K_{G_{\alpha \oplus \beta}}(X) \longrightarrow K_{G_{\alpha \oplus \beta}}(X)$$
will restrict to the direct summands to give a product
$$K_{G_{\alpha \oplus \beta},V}(X) \bigotimes K_{G_{\alpha \oplus \beta},W}(X) \longrightarrow K_{G_{\alpha \oplus \beta},V \otimes W}(X)$$
where $V$ and $W$ are irreducible complex $\mZ/n\mZ \times \mZ/m\mZ$-modules. Therefore, for a finite $G$-CW complex $X$,  the twisted product is a restriction of the standard product for the $G_{\alpha \oplus \beta}$-equivariant $K$-theory of $X$. The graded commutativity comes from the fact that $\alpha+\beta = \beta+\alpha$ as cocycles and the graded commutativity of the product in \cite{oliver-completion-theorem-2001}. 

\section{Twisted Bredon Cohomology and the Spectral Sequence}
This section will be devoted to the construction of a coefficient system for Bredon cohomology on the category of finite proper $G$-CW complexes. This gives \emph{twisted} Bredon cohomology which is related to twisted equivariant $K$-theory by a spectral sequence. 

\subsection{Twisted Bredon Cohomology}

Bredon cohomology was first introduced by Bredon \cite{bredon;equivariant-cohomology-theories} as a way to formulate equivariant obstruction theory.  Let $\mathcal{O}_G$ be the \emph{orbit category} of $G$; a category with one object $G/H$ for each subgroup $H \subseteq G$ and morphisms all $G$-maps between any two objects. Notice that there are morphisms $\phi:G/H \rightarrow G/K$ if and only if $H$ is subconjugate to $K$, since if $\phi(eH) = gK$, then $g^{-1}Hg \subseteq K$.

A \emph{coefficient system} for Bredon cohomology is a contravariant functor $M:\mathcal{O}_G^{op} \rightarrow \mathcal{A}b$ into the category of abelian groups. Given a $G$-CW complex $X$, there are some standard coefficient systems given by
$$\underline{C}_n(X)(G/H) = H_n( (X^{(n)})^H, (X^{(n-1)})^H; \mZ)$$
The connecting homomorphisms of the triple $((X^{(n)})^H, (X^{(n-1)})^H, (X^{(n-2)})^H)$ give rise to a boundary map
$$ \underline{ \partial}: \underline{C}_n(X) \rightarrow \underline{C}_{n-1}(X)$$
with $\underline{\partial}^2 = 0$. For two coefficient systems $M$ and $M'$, we denote by $Hom_{\mathcal{C}_G}(M, M')$ the abelian group of natural transformations between the two systems. This is an abelian group since the category $\mathcal{C}_G$ of coefficient systems is abelian \cite{may;equivariant-homotopy-cohomology}.
 
Now for any coefficient system $M$ and any $G$-CW compex $X$, we can define a cochain complex of abelian groups $C^*_G(X;M)$ as follows
$$C^n_G(X;M) = Hom_{\mathcal{C}_G}(\underline{C}_n(X), M), \mbox{ with } \delta=Hom_{\mathcal{C}_G}(\underline{\partial}, id).$$ 
This makes $C^*_G(X;M)$ a cochain complex of abelian groups. The homology of this complex is the \emph{Bredon cohomology} of $X$ with coefficient system $M$, denoted $H^*_G(X;M)$.

   Now, let $G$ be an arbitrary discrete group and $\alpha \in Z^2(G, U(1))$ be a normalized torsion cocycle. As before, to avoid cumbersome notation we will use $\alpha$ to denote both the cocycle and its restriction to any subgroup $H \subseteq G$ if it is clear from the context which is meant. We can define a coefficient system for Bredon cohomology of proper, finite $G-CW$ complexes by using the twisted  representation group functor $\mathcal{R_\alpha}(-)$. Define $\mathcal{R_\alpha}$ on objects by $\mathcal{R_\alpha}(G/H) = R_\alpha(H)$. Note that this is well defined if $H$ is a finite subgroup of $G$. Since we are considering proper action of discrete groups, the only orbits we are concerned with are of the form $G/H$ with $H$ finite. Therefore, this coefficient system is well defined in our case.

     The coefficient system is defined on the morphisms of $\mathcal{O}_G$ as 
follows. If $H \subseteq K$ are finite subgroups of $G$, then there is an 
obvious map $\phi: G/H \rightarrow G/K$ in the orbit category. We define 
$\mathcal{R}_\alpha(\phi)$ to be the map $R_{\alpha}(K) \rightarrow 
R_{\alpha}$ induced by restricting every representation to the 
subgroup. A map in the orbit category induced by $gHg^{-1} \subseteq K$ for 
some $g \in G$ is given by $g'H \rightarrow g'gK$, and it induces a map 
$R_{\alpha}(K) \rightarrow R_{\alpha}(H)$ by restricting every 
representation to the subgroup $gHg^{-1}$ and then using the isomorphism 
induced by conjugation. 

\begin{Def} The $\alpha$-twisted Bredon cohomology of a finite, proper 
$G-CW$ complex $X$ is defined to be $H^*_G(X; \mathcal{R_\alpha})$ where $\alpha \in Z^2(G, U(1))$ is a cocycle and $\mathcal{R_\alpha}$ 
is the coefficient system described above.
\end{Def}

It is worth mentioning that the name $\alpha$-twisted Bredon cohomology is a bit misleading. This is {\bf not} Bredon cohomology with twisted (local) coefficients. There is a theory of Bredon cohomology with local coefficients which generalizes the nonequivariant case \cite{mukherjee-pandey}. We use the name $\alpha$-twisted Bredon cohomology to emphasize the connection with $\alpha$-twisted equivariant $K$-theory which will be explored more deeply later.

Many common coefficient systems for Bredon cohomology extend to a quotient of $\mathcal{O}_G$ denoted by $Sub(G)$. $Sub(G)$ is the category of subgroups of $G$ and monomorphisms induced by subconjugation. $Sub(G)$ is a subquotient of $\mathcal{O}_G$ because all of the morphisms $G/H \rightarrow G/H$ coming from conjugation by an element of $C_G(H)$ are collapsed to the identity. So for a coefficient system to extend to $Sub(G)$, it must take any map $G/H \rightarrow G/H$ induced by conjugation by $g \in C_G(H)$ to the identity. It is clear that the coefficient system $\mathcal{R_\alpha}$ extends to $Sub(G)$, and so it defines a contravariant functor from $Sub(G)$ to $\mathcal{A}b$.
 
For the following important result we need to use induced modules.

\begin{Thm} \label{mackey}
The coefficient system $\mathcal{R_\alpha}$ has the structure of a Mackey functor to the category of $\mZ$-modules.
\end{Thm}

Before proving this, recall the definition of a Mackey functor. We will use the notation of \cite{luck}. Let $M:FGIN\mathcal{J}_G \rightarrow R-MOD$ be a \emph{bifunctor}. This means it is a pair
$(M_*, M^*)$ consisting of a covariant functor $M_*$ and a contravariant functor $M^*$ from $FGIN\mathcal{J}_G$ to $R-MOD$ which agree on objects. Here $FGIN\mathcal{J}_G$ is the category of finite subgroups of $G$ and injective homomorphisms induced by (sub)conjugation. For an injective group homomorphism $f:H \rightarrow K$ we will denote the map $M_*(f)$ by $ind_f$ and the map $M^*(f)$ by $res_f$. If $f$ is an inclusion of groups we will write $ind^{K}_H = ind_f$ and $res^H_K = res_f$. Such a bifunctor $M$ is a \emph{Mackey functor} if
\begin{itemize}
\item[(1)] For an inner automorphism $c(g):H \rightarrow H$, we have $M_*(c(g)) = id:M(H) \rightarrow M(H)$;
\item[(2)] For an isomorphism of groups $f:K \stackrel{\cong}{\longrightarrow} H$, the composites $res_f \circ ind_f$ and $ind_f \circ res_f$ are the identity;
\item[(3)] Double coset formula \\
           We have for two subgroups $H, K \subseteq L$
$$res^K_L \circ ind^L_H = \sum_{KgH \in K \backslash L/H} ind_{c(g):H \cap g^{-1}Kg \rightarrow K} \circ res^{H \cap g^{-1}Kg}_H$$
where $c(g)$ is conjugation with $g$.
\end{itemize}
\noindent \emph{Proof of Theorem \ref{mackey}.}
First, recall that when considering $\alpha$-twisted representations of a finite group $G$ it suffices to consider $\mC^{\alpha}G$-modules. We will use the notation of \cite{luck} as above.  We make $\mathcal{R_\alpha}$ into a bifunctor in the following way. Let its value on objects be the obvious one, $\mathcal{R_\alpha}(H) = R_\alpha(H)$. For the contravariant functor, let $(\mathcal{R_\alpha})^*(H \rightarrow K) = R_\alpha(K) \rightarrow R_\alpha(H)$ be the usual restriciton induced by inclusion (and possibly conjugation). For the covariant functor, let $(\mathcal{R_\alpha})_*(H \rightarrow K)$ be the map taking the $\mC^{\alpha}H$ -module $V$ to the $\mC^{\alpha}K$-module $\mC^{\alpha}K \otimes_{\mC^{\alpha}H} V^{(g)}$ where $gHg^{-1} \subseteq K$. Now the first two properties follow from the definitions, and the third follows from the Subgoup Theorem (Theorem \ref{sub}). So $\mathcal{R_\alpha}$ is a Mackey functor.
\hfill \qed
\vspace{0.2cm}

Since $R_\alpha(H)$ is not a ring, we cannot expect to have a multiplicative structure on $H^*_G(X, \mathcal{R_\alpha})$. However, we do have a sort of twisted multiplication
$$R_\alpha(H) \otimes R_\beta(H) \rightarrow R_{\alpha + \beta}(H)$$
as noted previously. This gives hope that there might be an analogous 
product structure on twisted Bredon cohomology, and in fact there is.

\begin{Thm} There is a pairing relating the $\alpha$-twisted and 
$\beta$-twisted Bredon cohomologies of a proper $G$-CW complex to its 
($\alpha + \beta$)-twisted Bredon cohomology. The pairing takes the form
$$\cup : H^n_G(X,\mathcal{R_\alpha}) \otimes H^m_G(X,\mathcal{R_\beta}) 
\rightarrow H^{n+m}_G(X,\mathcal{R_{\alpha + \beta}})$$
for $\alpha, \beta \in Z^2(G, U(1))$. This pairing is natural in $X$.
\end{Thm}

The proof of this theorem is exactly analagous to the proof of the existence of the cup product in singular cohomology. Also, since $\alpha + \beta = \beta + \alpha$ as cocycles, it is easy to show that the twisted product is graded commutative, i.e.
$$s \cup t = (-1)^{(nm)}t \cup s  \,\, \mbox{ for } s \in H^n_G(X, \mathcal{R_\alpha}), t \in H^m_G(X,\mathcal{R_\beta}).$$

It is clear that if $\alpha$ is the trivial cocycle, we have Bredon cohomology with coefficients in the representation ring functor. For simplicity call it untwisted Bredon cohomology. Notice that the product just defined makes $\alpha$-twisted Bredon cohomology a module over untwisted Bredon cohomology in a natural way. This is completely analagous to what happened with $\alpha$-twisted equivariant $K$-theory. This isn't suprising, since we can connect the two theories by a twisted equivariant Atiyah Hirzebruch spectral sequence.

\subsection{The Spectral Sequence}
Next we see that there is a spectral sequence connecting the $\alpha$-twisted Bredon cohomology and the $\alpha$-twisted equivariant $K$-theory of finite proper $G$-CW complexes. This spectral sequence is a special case of the more general spectral sequence constructed by Davis and L\"{u}ck \cite{davis-l*;spaces-category-spaces-assembly-isomorphism-conjectures}.

\begin{Thm}
Let $X$ be a finite proper $G$-CW complex for a discrete group $G$,
and let $\alpha \in Z^2(G, U(1))$ be a normalized torsion cocycle.
Then there is a spectral sequence  with 
$$E_2^{p,q} = \left\{ \begin{array}{ll} H^p_G(X,\mathcal{R}_\alpha) & \mbox{ if } q \mbox{ is even} \\
0 & \mbox{ if } q  \mbox{ is odd} \\
\end{array} \right.$$
so that $E^{p,q}_{\infty} \Longrightarrow
\,^{\alpha}K_G^{p+q}(X)$.
\end{Thm}

The construction of the spectral sequence is completely analogous to the construction of the standard Atiyah-Hirzebruch spectral sequence \cite{atiyah-hirzebruch;vector-bundles-homogeneous-spaces} using the skeletal filtration of $X$. In fact, this spectral sequence is a direct summand of the spectral sequence in \cite{l*-oliver;chern-proper-equivariant-characters}. We can see this using the fact that $R_\alpha(G)$ is a direct summand of $R(G_\alpha)$ and $\,^\alpha K^*_G(X)$ is a direct summand of $K^*_{G_\alpha}(X)$ because of the splitting (\ref{splitting}). The maps in the spectral sequence considered in \cite{l*-oliver;chern-proper-equivariant-characters} clearly respect the splitting since they are induced by the skeletal filtration, and so our spectral sequence comes out as a direct summand.

Note that the $E_2$ page of the spectral sequence for $\alpha$-twisted equivariant $K$-theory can be seen as a module over the $E_2$ page of the spectral sequence for untwisted equivariant $K$-theory using the product described in Section 4.1. Also, it has been shown \cite{l*-oliver;chern-proper-equivariant-characters} that the spectral sequence for untwisted equivariant $K$-theory is a multiplicative spectral sequence. This can be used to make the spectral sequence for $\alpha$-twisted $K$-theory into a module over the spectral sequence for untwisted equivariant $K$-theory since all of the differentials respect the module structure.

We will use a result of L\"{u}ck to show that this spectral sequence collapses modulo torsion. In \cite{luck}, he constructs a Chern character for any proper equivariant cohomology theory which maps from the theory to its associated Bredon cohomology. In our case, $\alpha$-twisted equivariant $K$-theory is the proper equivariant cohomology theory and its associated Bredon cohomology is $\alpha$-twisted Bredon cohomology. In fact, since our spectral sequence is a summand of the spectral sequence from \cite{l*-oliver;chern-proper-equivariant-characters} for $G_\alpha$-equivariant $K$-theory of $X$, the Chern character for $\alpha$-twisted $K$-theory can be seen as a restriction of the Chern character for $G_\alpha$-equivariant $K$-theory.

To apply \cite[Theorem 5.5]{luck} which gives that the Chern character is an isomorphism, it suffices to show that $\mathcal{R}_\alpha(-) \otimes \mQ$ has a Mackey structure. By Theorem \ref{mackey}, $R_\alpha(-)$ has a Mackey structure as a functor into $\mZ-MOD$. Since $\mQ$ is a flat $\mZ$-module, tensoring with $\mQ$ preserves exact sequences; in particular it preserves isomorphisms. This means that the three conditions of a Mackey functor are still preserved. So $R_\alpha(-) \otimes \mQ$ does have a Mackey structure.

After tensoring with $\mQ$, the equivariant Chern character gives an isomorphism between the $E_2$ page and the $E_\infty$ page of the spectral sequence. This tells us that the spectral sequence collapses at the $E_2$ page after tensoring with $\mQ$.

\section{Connections to other constructions}
There has been a lot of activity in studying twisted equivariant $K$-theory recently motivated by orbiforlds \cite{adem-ruan}, \cite{bunke-schick;T-duality}, \cite{tu-xu-lg;differentiable-stacks}. Generally, as in \cite{atiyah-segal;twisted}, if $X$ is a $G$-space, the twistings considered are classified by $H^3(X \times_G EG;\mZ)$. In fact, it is a little more subtle than that because the twistings actually depend on cocycles not cohomology classes. Two cohomologus cocycles do give isomorphic twisted equivariant $K$-groups. However the isomorphism is not canonical; it depends upon the coboundary connecting the two cocycles. 

The twistings discussed here fit into the more general setting in a very natural way. Given a discrete group $G$ and a normalized torsion cocycle $\alpha \in Z^2(G,U(1))$, $\alpha$ defines a cohomology class $[\alpha] \in H^2(G,U(1))$. The short exact sequence of coefficients
$$ 1 \rightarrow \mZ \rightarrow \mR \rightarrow U(1) \rightarrow 1$$ 
gives rise to a long exact sequence in the cohomology of $G$. In particular, it gives a connecting homomorphism
$$
H^2(G,U(1)) \stackrel{\delta}{\longrightarrow} H^3(G,\mZ).
$$
Now given a finite $G$-CW complex $X$, we have the fibration
$$\xymatrix@1{
X \ar[r] & X \times_G EG \ar[r]^p & BG \\
}$$
which induces a map $p^*:H^3(BG;\mZ) \rightarrow H^3(X \times_G EG; \mZ)$. So twistings of the type we have discussed are a subclass of the more general twistings. 

Morally, the $G$-equivariant $K$-theory of a space $X$ is a combination of the K-theory of the space, $K(X)$, and the complex representation ring of the group, $R(G)$. General twistings involve both pieces, however the twistings we have considered only depend on the group. They can be interpreted as leaving $K(X)$ alone and only twisting the ingredient $R(G)$.

Another approach to twisted $K$-theory for orbifolds come from the groupoid description of an orbifold. In \cite{tu-xu-lg;differentiable-stacks}, given an orbifold groupoid $\mathfrak{X}$, the twisted orbifold $K$-theory of $\mathfrak{X}$ is described in terms of the $K$-theory of a particular $C^*$-algebra. In the case where the orbifold is a quotient orbifold given by a proper action of a Lie Group $G$ on a $G$-compact space $X$ and the twisting is trivial, this resticts to the $G$-equivariant $K$-theory of $X$. For quotient orbifolds the twistings again are classified by $H^3(X \times_G EG; \mZ)$, and twisted orbifold $K$-theory is just the twisted equivariant $K$-theory of $X$. One important question raised in \cite{tu-xu-lg;differentiable-stacks} is : When can the twisted orbifold $K$-theory be realized in terms of finite rank bundles? In \cite{tu-xu-lg;differentiable-stacks}, they prove that if the twisting is not torsion, then it is impossible to describe twisted orbifold $K$-theory just using finite rank bundles. However, if the twisting is torsion, it is unknown if it is possible or not.

They prove \cite[Theorem 5.28]{tu-xu-lg;differentiable-stacks} a general result in this direction. For a compact quotient orbifold and a torsion twisting, the theorem says that if there is a single finite rank twisted bundle, then the twisted orbifold $K$-theory is representable in terms of finite rank bundles. They conjecture that for any torsion twisting, there will always exist at least one finite rank twisted bundle. So our construction gives a partial result for their conjecture. If $G$ is a discrete group, $X$ a finite proper $G$-CW complex, then given any torsion twisting $[\alpha] \in H^3(X \times_G EG; \mZ)$ with $[\alpha]$ in the image of $p^* \circ \delta$ there is a finite rank $\alpha$-twisted bundle (Theorem \ref{bundle}).

Tu and Xu have also constructed a Chern character for twisted orbifold $K$-theory which lands in a twisted version of deRahm cohomology \cite{tu-xu;chern-orbifolds}. Since their definition of twisted $K$-theory involves operator algebras, their Chern character factors through the periodic cyclic homology groups of the operator algebra. In the case of a discrete group acting properly and a torsion twisting coming from the group, the Chern character mentioned in this paper gives an alternate method of computation for the twisted $K$-theory groups.

\end{document}